\documentclass[]{interact}

\usepackage[T1]{fontenc} 
\usepackage[utf8]{inputenc} 
\usepackage{floatrow} 
\usepackage{tablefootnote}
\usepackage[linktoc=page]{hyperref}
\hypersetup{linkcolor  = black}
\hypersetup{citecolor=black} 
\hypersetup{colorlinks=true}
\hypersetup{urlcolor=cyan}
\usepackage{color, colortbl}
\usepackage{graphicx, amsmath, amsthm, amssymb, mathrsfs, amscd}
\usepackage{tablefootnote}
\usepackage{footnote}
\usepackage{epstopdf}
\usepackage[caption=false]{subfig}
\usepackage[numbers,sort&compress]{natbib}
\bibpunct[, ]{[}{]}{,}{n}{,}{,}
\makeatletter
\def\NAT@def@citea{\def\@citea{\NAT@separator}}
\makeatother
\theoremstyle{plain}
\newtheorem{theorem}{Theorem}[section]
\newtheorem{lemma}[theorem]{Lemma}

\theoremstyle{definition}
\newtheorem{definition}[theorem]{Definition}

\theoremstyle{remark}
\newtheorem{remark}{Remark}

\usepackage{lineno,hyperref}
\modulolinenumbers[5]
\usepackage{tablefootnote}
\hypersetup{linkcolor  = blue}
\hypersetup{citecolor=blue} 
\hypersetup{colorlinks=true}
\hypersetup{urlcolor=cyan}
\usepackage{color, colortbl}
\usepackage{multirow}
\usepackage{graphicx, amsmath, amsthm, amssymb, mathrsfs, amscd}
\usepackage{caption}
\usepackage{enumerate}
\usepackage{algorithm}
\usepackage[noend]{algcompatible}
\usepackage{url}
\usepackage{tikz}
\usepackage[toc,page]{appendix}
\usepackage{footnote}

\usepackage{color}
\usepackage{color, colortbl}
\usepackage{tablefootnote}
\usepackage{footnote}
\theoremstyle{definition}

\theoremstyle{definition}

\theoremstyle{definition}

\theoremstyle{definition}
\theoremstyle{definition}
\makeatletter
\def\cleardoublepage{\clearpage\if@twoside \ifodd\c@page\else
  \hbox{}
  \vspace*{\fill}
    \vspace{\fill}
  \if@twocolumn\hbox{}\newpage\fi\fi\fi}
\makeatother

\begin{document}
\title{The Image of Critical Circle and Zero-free Curve for Quadrinomials}
\author{
\bigskip \name{Oluma~Ararso~Alemu \textsuperscript{1,\textcolor{blue}{*}}, and Hunduma~Legesse~Geleta \textsuperscript{1}}\footnote{ \textcolor{blue}{*}Corresponding author}
 \affil{\textsuperscript{1} Department of Mathematics, Addis Ababa University, Addis Ababa, Ethiopia. \\  Email: oluma.ararso@aau.edu.et and hunduma.legesse@aau.edu.et }} 
\maketitle
\begin{abstract}
The location of the zeros of a two-parameter family of complex-valued harmonic quadrinomials depends on the parameters.  In this paper, we determine and demonstrate that  the image of some critical circle under these two-parameter family of complex-valued harmonic quadrinomial is a hypocycloid. We also determine a zero free curve for two-parameter family of  quadratic quadrinomial. 
\end{abstract}
\section{\textcolor{blue}{Introduction}}
An open connected set $G$ is simply connected if and only if for each harmonic function $u$ on $G$ there is a harmonic function $v$ on $G$ such that the function $f(x,y) = u(x,y) + iv(x,y)$ is analytic on $G.$ This implies that every harmonic function on a simply connected region has a harmonic conjugate.
A continuous function $f(x,y)=u(x,y)+iv(x,y)$ defined in the region $G$ is said to be a complex valued harmonic function  if $u$ and $v$ are real valued harmonic (but not necessarily conjugates to each other)  in $G.$ Note that one way of thinking of a function $f(x,y) = u(x,y) + iv(x,y)$ as being analytic is that $f$ can be expressed as a function of $z=x+iy$ only without using $\overline{z} = x-iy.$ For instance, the function $f(z)= z^2$ is analytic while $f(z)= z\overline{z}$ is not. As a result, a one-to-one mapping $f:G \rightarrow \Omega, f(z) = u(x,y) + iv(x,y)$ from a region $G$ in $xy- plane$ to a region $\Omega$ in $uv-plane$ is harmonic mapping if the two coordinate functions $u(x,y)$ and $v(x,y)$ are real harmonic(but not necessarily \bigskip conjugates) in $G$. \\
A useful result in determining the location of the zeros of analytic polynomials is the argument principle. Recall that if $f$ is analytic inside and on positively oriented rectifiable Jordan curve C and  $f(z) \neq 0$  on C, then  the winding number of the image curve $f(C)$ about the origin $,\frac{1}{2\pi}\Delta _C~ \mathrm{arg} f(z),$  equals the total number of zeros of $f$ in D, counted according to multiplicity where D is a plane domain bounded by C. In 1940 Kennedy  \cite{kennedy1940bounds} showed that  the roots of analytic trinoimial equations of the form $z^n + az^k + b = 0 ,$ where $ab \neq 0,$  have  certain bounds to their respective absolute values. In 2006, Dehmer \cite{dehmer2006location} proved that all the zeros of complex-valued analytic polynomials lie in certain closed disk.
 In 2012,  Melman \cite{melman2012geometry} investigated the regions in which the zeros of analytic trinomials of the form $p(z) = z^n - \alpha z^k - 1,$ with integers $n \geq 3$ and $1 \leq k \leq n - 1$ with $\mathrm{gcd}(k, n) = 1,$ and $\alpha \in \mathbb{C},$ lie. He determined the zero inclusion regions for the following cases: a$)$ for any value of $|\alpha|;$ $b)$  $| \alpha| > \sigma (n, k);$ and c$)$  $| \alpha| < \sigma (n, k),$ where  $\sigma(n, k)=\frac{n}{n-k} (\frac{n-k}{k})^{\frac{k}{n}}.$ In each cases, Melman provided useful information  on \bigskip the  location  of the zeros of $p.$ \\
One area of investigation that has recently become of interest is the number and location of zeros of complex-valued harmonic polynomials. If $f(x,y) = u(x,y) + iv(x,y)$ is harmonic in simply connected domain $G,$ then we can decompose $f$ as $f(z)=h(z)+\overline{g(z)}$  where $g$ and $h$ are a complex valued analytic functions in $G.$ In the decomposition of a complex valued harmonic polynomial $f(z)= h(z) +\overline{g(z)},$ $h(z)$ is called analytic part and $g(z)$ is said to be co-analytic part of $f.$ The function $f$ is said to have degree $n$ if either $a_n$ or $a_{-n}$ is non-zero.The Jacobian of complex valued function $f(z) = u+iv$ is defined by
 $ J_f =
 \begin{vmatrix}
 u_x & v_x\\ u_y & v_y
\end{vmatrix} = u_xv_y-u_yv_x. 
$
For analytic function $f=u+iv, J_f = |f'|^2$ and for harmonic function $f=h+\overline{g}, J_f = |h'|^2-|g'|^2.$
A complex valued harmonic function $f(z)=h(z)+\overline{g(z)}$ is said to be sense preserving at $z$ if $J_f(z) >0$ and is sense reversing at $z$ if $J_f(z) < 0.$ If $f$ is neither sense preserving nor sense reversing at $z,$ then $f$ is said to be singular at $z.$ All non-constant analytic function is a sense-preserving. The dilatation of a complex valued harmonic function $f(z)=h(z)+\overline{g(z)}$ is defined to be  
$\omega(z) = \frac{g'(z)}{h'(z)}.$
Note that the dilatation of analytic function is zero. For generality, the dilatation is a measure of how far a harmonic function is from being analytic. Next we generalize that the Jacobian of complex valued harmonic function cannot be zero whenever the function is locally univalent, which was proved by Hans Lewi in \cite{lewy1936non}. \\

It was shown by Bshouty \textit{et al.} \cite{bshouty1995exact} that there exists  a complex-valued harmonic polynomial $ f = h + \overline{g},$ such that $h$  is an analytic polynomial of degree $n,$ $g$  is an analytic polynomial of degree $m < n$ and $f$ has exactly $n^2$ zeros counting with  multiplicities in the field of complex numbers,  $\mathbb{C}.$ We are motivated by the work of Brilleslyper et al. \cite{brilleslyper2020zeros},  Kennedy \cite{kennedy1940bounds}, and Dehmer \cite{dehmer2006location} on the number and location of zeros of  trinomials. In 2020, Brilleslyper et al. \cite{brilleslyper2020zeros} studied the number of zeros of harmonic trinomials of the form $p_c(z) = z^n + c\overline{z}^k - 1$  where $1 \leq k \leq n - 1, n \geq 3, c \in \mathbb{R^+},$ and $\mathrm{gcd}(n, k) = 1.$ They showed that the number of zeros of  $p_c(z)$ changes as $c$ varies and proved that the distinct number of zeros of $p_c(z)$ ranges from $n$ to $n+2k$. Among other things, they used the argument principle for harmonic function that can be formulated as a direct generalization of the classical result for analytic functions (Duren et al. \cite{duren1996argument}). Recently,   Hunduma Legesse and Oluma Ararso \cite{legesse2022location} solved an interesting problem stated in \cite{brilleslyper2020zeros} by Brilleslyper et al. We found the zero inclusion region of complex-valued harmonic trinomials. Not only finding the zero inclusion regions of these particular trinomials, the region containing the zeros of any complex-valued harmonic polynomials is discovered in the same paper. Moreover, the same authors in \cite{alemu2022zeros} studied a special type of complex-valued harmonic polynomial of the form $Q_{b,c}(z)= bz^k+\overline{z}^n+c\overline{z}^m+z$  with $b,c \in \mathbb{R}\setminus\{0\},$ $k\geq n>m,$ and $k,n,m \in \mathbb{N}$. They have determined the curve separating sense-preserving region from sense-reversing region. This curve is defined as $\Upsilon _{b,c}= \left\lbrace z \in \mathbb{C}:|z|=\mathcal{M}_{b,c} =\left(\frac{c^2-1}{k^2(b^2-1)}\right)^\frac{1}{2k-2} \right\rbrace$ and is sad to be a \textbf{critical circle.} An interesting open problems raised in this article was:
\begin{itemize}
\item[(1)] How many roots of $q(z)$ are $\mathcal{M}_{b,c}$\textit{-modular}? What can be said on the number of zeros inside and outside of the circle $\Gamma _{b,c}= \{z \in \mathbb{C}:|z|=\mathcal{M}_{b,c} \}?$
\item[(2)] What is image of the circle $\Gamma _{b,c}= \{z \in \mathbb{C}:|z|=\mathcal{M}_{b,c} \}$ under the quadrinomial $q(z)=bz^k+\overline{z}^n +c\overline{z}^m +z?$
\end{itemize}

In this paper, we first look at  a two-parameter complex-valued harmonic quadrinomials $Q_{b,c}(z)= bz^k+\overline{z}^n+c\overline{z}^m+z$  with $b,c \in \mathbb{R}^+\setminus\{1\},$ $k\geq n>m,$ and $k,n,m \in \mathbb{N}$ and find the image of quadrinomic critical circle under this quadrinomial using some considerations stated in \cite{alemu2022zeros}. We also find the number of some parametric radial zeros.  The main theorems in this paper   are  Theorem $\ref{r1}$ and Theorem $\ref{r2}.$  The structure of our paper is  as follows. In section $\ref{p}$, we present some important preliminary results that will be used to prove our two main theorems. In section $\ref{r}$, we state and prove the main results of this paper.   By using mathematical software,  Theorem $\ref{r1}$  demonstrates that the image of circle under $Q_{b,c}(z)$ is a hypocycloid with $k+1$-cusps. Theorem $\ref{r2}$ determines the quadrinomial $Q_{b,c}(z)= bz^k+\overline{z}^n+c\overline{z}^m+z$  with $b,c \in \mathbb{R}^+\setminus\{1\},$ $k\geq n>m,$ and $k,n,m \in \mathbb{N}$ has no zeros on the curve separating sense-preserving and sense-reversing region. In section $\ref{c},$ we come up with conclusion.
   \section{\textcolor{blue}{Preliminaries}}$\label{p}$
   In this section we review some important concepts and results that we will use later on to prove our theorems. We begin by stating the well known results about hypocycloid, \bigskip zeros of complex-valued harmonic polynomials and some useful theorems and lemmas.\\
   The curve in the plane that separates the sense preserving and sense reversing region for complex-valued harmonic polynomials plays a significant role in complex analysis. A straightforward computation gives that $|\omega (z)| =1$ if and only if $|z|= \left[\frac{(c^2-1)}{k^2(b^2-1)}\right]^{\frac{1}{2k-2}},$ with $Q_{b,c}(z) = Q_{b,c}(z)= bz^k+\overline{z}^k+c\overline{z}+z$ being sense-reversing on the interior of this circle and sense preserving on its exterior. This is proved in \cite{alemu2022zeros}.
\begin{definition}
A hypocycloid centered at the origin is the curve traced by a fixed point on a circle of radius $r$ rolling inside a larger origin-centered circle of radius $R.$ \\ The curve is given by the parametric equations 
$$ x(\phi ) = (R-r)cos(\phi )+rcos\left( \frac{R-r}{r}\phi \right);$$ $$y(\phi ) = (R-r)sin(\phi )-rsin\left( \frac{R-r}{r}\phi \right) .$$
\end{definition}
If the ratio $\frac{R}{r}$ is written in reduced form as $\frac{p}{q} \in \mathbb{Q},$ then the hypocycloid has $p~cusps$ and each arc connects cusps that are $q$ away from each other in a counter clockwise direction. Such a hypocycloid is called a $(p,q)$ hypocycloid, and  the range of $\phi$ values to trace the entire hypocycloid is $0 \leq \phi \leq 2\pi q.$
\begin{definition}
If both analytic and co-analytic part of a two-parameter family of quadrinomial $Q_{b,c}(z)= bz^k+\overline{z}^n+c\overline{z}^m+z,$ with $b,c \in \mathbb{R}^+\setminus\{1\},$ $k\geq n>m,$ and $k,n,m \in \mathbb{N},$ are quadratic polynomials, then $Q_{b,c}(z)$ is called \textbf{quadratic quadrinomial.}
\end{definition} 
The following theorem is proved in \cite{hayami2011hypocycloid}.
\begin{theorem}

Let $f(z)$ be holomorphic in the closed unit disk $\overline{\mathbb{U}},$ with $f'(z) \neq 0$ for
 $ z\in \overline{\mathbb{U}},$ and let $F(t)=(n+1)t+2\mathrm{arg}(f'(e^{it}))~~~ (- \pi \leq t < \pi), ~~~ n=2,3,4,\dot{...}.$ If for each $k \in K = \left\{ 0, \pm 1, \pm 2,..., \pm \left[ \frac{n+3}{2} \right]\right\}$ where $[~~]$ denotes the Gauss symbol, the equation $F(t) = 2k\pi$ has at most a single root in $[-\pi, \pi)$ and for any $k \in K$ there exist exactly $(n + 1)$ such roots in $[-\pi, \pi),$ then the harmonic function $f(z) = h(z) + \overline{g(z)}$ with $g'(z)=z^{n-1}f'(z)$ is univalent in $\overline{\mathbb{U}}$, sense preserving and the image of $\overline{\mathbb{U}}$ by $h(z)$ is a hypocycloid of $n + 1$ cusps.
\end{theorem}

\begin{definition}
A one to one complex-valued harmonic function is said to be univalent harmonic function. A \textit{locally univalent} functions are functions that are one to one locally. That means, a function $f(z)$ is said to be locally univalent in the domain $G$ if there is a small neighborhood around each point $z_0 \in G,$ such as a smaller disk centered at $z_0$ and the function is one to one in that neighborhood.
\end{definition}
\begin{definition}
A complex-valued harmonic function $f(z)=h(z)+\overline{g(z)}$ is said to be \textit{sense-preserving} at $z_0$ if $J_f(z_0) >0$ and is \textit{sense-reversing} at $z_0$ if $J_f(z_0) < 0,$ where $J_f(z_0)$ is the jacobian of $f$ which is given by $J_f =
 \begin{vmatrix}
 u_x & v_x\\ u_y & v_y
\end{vmatrix} = u_xv_y-u_yv_x.$ If $f$ is neither sense-preserving nor sense-reversing at $z_0,$ then $f$ is said to be \textit{singular} at $z_0.$
\end{definition}
\begin{definition}
The \textit{dilatation} of a complex-valued harmonic function $f(z)=h(z)+\overline{g(z)}$ is defined to be  
$
 \omega(z) = \frac{g'(z)}{h'(z)}. 
$
\end{definition}
\begin{theorem}\textbf{(\cite[Lewy's Theorem] {lewy1936non})}\label{2.1} 
If $f$ is a complex-valued harmonic function that is locally univalent in a domain $\mathbb{D} \subset \mathbb{C},$ then its Jacobian $,J_f(z),$ never vanish for all $z \in \mathbb{D}.$ 
\end{theorem}
As an immediate consequence of Theorem $\ref{2.1}$, a complex-valued harmonic function $f(z)=h(z)+\overline{g(z)}$ is locally univalent and sense-preserving if and only if $h'(z) \neq 0$ and $|\omega (z)|<1.$ A harmonic function $f(z) = h(z) + \overline{g(z)},$ is called sense-preserving at $z_0$ if the Jacobian $J_f(z) > 0$ for every $z$ in some punctured neighborhood of $z_0.$ We \bigskip also say that $f$ is sense-reversing if $\overline{f}$ is sense-preserving at $z_0.$ \\
 Recall to Argument principle for harmonic functions in \cite{duren1996argument}: Let $f$ be a complex-valued harmonic function in a Jordan domain $\mathbb{D}$ with boundary $\Gamma .$ Suppose that $f$ is continuous in $\overline{\mathbb{D}}$ and $f \neq 0$ on $\Gamma .$ Suppose that $f$ has no singular zeros in $\mathbb{D},$ and let $N = N_+ - N_- ,$ where $N_+$ and $N_-$ are the number zeros in sense-preserving region and sense-reversing region of $f$ in $\mathbb{D},$ respectively. Then, $\bigtriangleup _\Gamma \mathrm{arg} f(z)  = 2 \pi N.$
\begin{definition}
Suppose $z_0$ is a fixed point of a function $f,$ that is, $f(z_0)=z_0.$ The number $\lambda = |f'(z_0)|$ is called the \textit{multiplier} of $f$ at $z_0.$ We classify the fixed point according to $\lambda$ in the following definition.
\end{definition}
\begin{definition}
A point $\zeta$ is called a critical point of a polynomial $f$ if $f'(\zeta) = 0.$  A fixed point $\zeta \in \mathbb{C}$ is \textit{attractive, repelling} or \textit{neutral} if, respectively, $|f'(\zeta)| < 1, |f'(\zeta)| > 1$ or
$|f'(\zeta)| = 1.$  A neutral fixed point is rationally neutral if $f'(\zeta)$ is a root of unity. We shall say that a fixed point $\zeta$ attracts some point $w \in \mathbb{C}$ provided that the sequence $f^k(w) = \underbrace{f(w) \circ f(w) \circ \cdots \circ f(w)}_{k-copies}$ converges to $\zeta.$ 
\end{definition}
\begin{definition}
Let $f(z)$ be a polynomial   such that the conditions $ |f'(z_0)| = 1$ and $\overline{f(z_0)} = z_0$ are not satisfied simultaneously for any $z_0 \in \mathbb{C}.$ Then $f(z)$ is said to be a \textit{regular} polynomial.
\end{definition}
\begin{definition}\cite{alemu2022zeros}
The roots of the quadrinomial $q(z)=bz^k+\overline{z}^n +c\overline{z}^m +z$ that lie on the circle of radius $\mathcal{M}_{b,c}$ are said to be the $\mathcal{M}_{b,c}$\textit{-modular} roots.
\end{definition}
\begin{definition}
The grand orbit  of a point $z \in \hat{\mathbb{C}},$ denoted by $[z],$ under a map $f$ is defined as the set of all points $w$ such that there exists $m, n \geq ~0$ with $f^m(z) = f^n(w).$
\end{definition}
\begin{remark} A grand orbit under a transformation $T$ is an equivalence class of the relation $x \sim y~~ iff ~~T^p(x) = T^q(y)$ for some $p, q > 0.$
\end{remark}
   Note that according to Descartes's Rule of Signs, polynomial $p(x)$ has no more positive roots than it has sign changes and  has no more negative roots than  $p(-x)$  has sign changes. A polynomial may not achieve the maximum allowable number of roots given by the Fundamental Theorem, and likewise it may not achieve the maximum allowable number of positive roots given by the Rule of Signs. The upper bound can always be attained for any sign sequence of coefficients that contains no zeros, and an explicit formula for finding such a polynomial can easily be found and this polynomial can be modified to reduce the number of roots by any even number while maintaining the same sign sequence. 
\begin{theorem}\textbf{ (\cite [Bezout's Theorem In the Plane]{kirwan1992complex})} $\label{II0}$
  Let $f$ and $g$ be relatively prime polynomials in the real variables $x$ and $y$ with real coefficients, and let $\mathrm{deg}h = n$ and $\mathrm{deg}g = m.$ Then the two algebraic curves $f(x, y) = 0$ and $g(x, y) = 0$ have at most $mn$ points in common.
  \end{theorem}
  Bezout's theorem is one of the most fundamental results about the degrees of polynomial surfaces and it bounds the size of the intersection of polynomial surfaces.
\begin{theorem}\textbf{(\cite[Wilmshurst]{wilmshurst1998valence})}$\label{II1}$
If $f(z)= h(z)+\overline{g(z)} $ is a harmonic polynomial such that $\mathrm{deg}h =n>m=\mathrm{deg}g$ and $\lim_{z \rightarrow \infty}f(z) = \infty,$ then $f(z)$ has at most $n^2$ zeros.
\end{theorem}
The proof of Theorem $\ref{II1} $ can also readily follows from Bezout's theorem. Let $f(z)=h(z) +\overline{g(z)}= 0$  with $h=u+iv$ and $g=\alpha +i\beta.$ Then $f(z)= \left( u(x,y)+\alpha(x,y)\right)+i\left(v(x,y)-\beta(x,y)\right)=0+i0.$ Now define $\eta (x,y)
:= u(x,y)+\alpha(x,y)$ and $\zeta (x,y):=v(x,y)-\beta(x,y).$ Then we do have a homogeneous polynomial equations equation 
$$
\begin{cases}
\eta (x,y)=0\\
 \zeta (x,y)=0.
 \end{cases}
$$
 Here we know that $\mathrm{deg} \eta =n$ and $\mathrm{deg}\zeta=n.$ Therefore by Bezout's Theorem the \bigskip maximum number of zeros of $\eta$ and $\zeta$ in common is $n^2.$ 
\begin{theorem}\textbf{(\cite[Dehmer, M.]{dehmer2006location})}$\label{II3}$
   Let $f(z) = a_nz^n+a_{n-1}z^{n-1}+ \cdots +a_1z+a_0, ~~~ a_n \neq 0$ be a complex polynomial. All zeros of $f(z)$ lie in the closed disc $K(0,\mathrm{max}(1,\delta)),$ where $M:= \mathrm{max}\{ |\frac{a_j}{a_n}| \} ~~\forall j=0,1,2,\cdots ,n-1$ and $\delta \neq 1$ denotes the positive root of the equation $z^{n+1} - (1+M)z^n +M =0.$ 
   \end{theorem}
   Note that this theorem is also another classical result for the location of the zeros of analytic complex polynomial which depends on algebraic equation's positive root. Descartes' Rule of Signs plays a great role in the proof as shown by Dehmer. 
   \begin{lemma}$\label{II4}$
  A complex valued harmonic function $f(z)=h(z)+\overline{g(z)}$ is analytic if and only if $f_{\overline{z}} =0.$
   \end{lemma}
   
    The following theorem  is stated and proved in \cite{legesse2022location}.
   \begin{theorem}$\label{a'}$
 Let $$p(z) = h(z) + \overline{g(z)} = a_0+ \sum_{k=1}^n a_kz^k + \overline{\left(b_0+\sum_{k=1}^m b_kz^k\right)}$$ be a complex-valued harmonic polynomial with $\mathrm{deg}h(z) = n \geq m=\mathrm{deg}g(z)$ and let  $r \neq 1$ be positive real root of equation $~~x^{n+1} -(1+M)x^n+M=0~~$ where $M:= \mathrm{max}\left\{ \frac{|a_j|+|b_j|}{|a_n|}\right\}_{j=0}^{n-1}.$  Then all the zeros of $p(z)$ lie in the closed disk $D(0,R)$ where $R = \mathrm{max}(1, r).$ 
 \end{theorem}
 Note that Theorem $\ref{a'}$ is true for every complex-valued harmonic polynomial and for our quadrinomial, the following theorem is stated and proved in \cite{alemu2022zeros}.
 
 \begin{theorem}$\label{a''}$
Let $q(z) = bz^k+\overline{z}^n +c\overline{z}^m +z$ where $k,m, n \in \mathbb{N}$ with $n>m$ and $b,c \in \mathbb{R}\setminus \{0\}.$  Let $|c|>1$ and $k>n.$ If $\delta _1 \neq 1$ is the positive real root of the equation
$$
|b|x^{k+1}-(|b|+|c|)x^k+|c| =0,
$$
 then all the zeros of $q(z)$ lie in the closed disk $D(0;R)$ where $R= \mathrm{max}\{ 1, \delta _1 \}.$
\end{theorem}
\section{\textcolor{blue}{Main Results}}$\label{r}$
Under this section, we find the image of critical circle under complex-valued harmonic quadrinomials and we find the number of zeros on the critical circle.
\begin{theorem}\label{r1}
Let $Q_{b,c}(z)= bz^k+\overline{z}^n+c\overline{z}^m+z$ be a two-parameter family of complex-valued harmonic quadrinomial with $b,c \in \mathbb{R}^+\setminus \{ 1\}, k\geq n > m$ and $k,m,n \in \mathbb{N}.$ The image of quadrinomic critical circle $$\Upsilon(z)=\left\{z:|z|=\mathcal{M}_{b,c} = \left[\frac{(c^2-1)}{k^2(b^2-1)}\right]^{\frac{1}{2k-2}} \right\}$$ under quadrinomial $Q_{b,c}(z) $ is a $\left(k+1, k \right)$-hypocycloid centered at $ 2\left(\frac{R-r)}{b+1}\right)sin\phi - \left(\frac{2r}{b+1}\right)sin\left( \frac{R-r}{R}\phi\right)$ where $b=c, k=n$ and $m=1.$ 
\end{theorem}
\begin{proof} Let $Q_{b,b}\left( \Upsilon_{b,b}\right)=\mathcal{X}(\phi)+i\mathcal{Y}(\phi).$
 Then 
$$ Q_{b,b}\left( \Upsilon_{b,b}\right)= Q_{b,b}\left( \mathcal{M}_{b,b}e^{i\theta}\right) = b\left(\mathcal{M}_{b,b}e^{i\theta}\right)^k+\left(\overline{\mathcal{M}_{b,b}e^{i\theta}}\right)^k+b \left(\overline{\mathcal{M}_{b,b}e^{i\theta}}\right) +\mathcal{M}_{b,b}e^{i\theta}~~~~~~~~~~~~~~~~~~$$ 

$$=b\left(\mathcal{M}_{b,b}\right)^k e^{ik\theta}+\left(\overline{\mathcal{M}_{b,b}}\right)^ke^{-ik\theta} + b\left(\overline{\mathcal{M}_{b,b}}\right)^ke^{-ik\theta} +\mathcal{M}_{b,b}e^{i\theta} ~~~~~~~~~~$$

 $$= b\left( \frac{1}{k^2}\right)^{\frac{k}{2k-2}}e^{ik\theta} + \left( \frac{1}{k^2}\right)^{\frac{k}{2k-2}}e^{-ik\theta} + b\left( \frac{1}{k^2}\right)^{\frac{1}{2k-2}}e^{-i\theta} + \left( \frac{1}{k^2}\right)^{\frac{1}{2k-2}}e^{i\theta}$$
 
 $$=b\left( \frac{1}{k^2}\right)^{\frac{k}{2k-2}}\left( cosk\theta + isink\theta \right)+ \left( \frac{1}{k^2}\right)^{\frac{k}{2k-2}}\left( cosk\theta - isink\theta \right) + ~~~~~~~~~~~~~$$ $$ b\left( \frac{1}{k^2}\right)^{\frac{1}{2k-2}}\left( cos\theta - isin\theta \right) + \left( \frac{1}{k^2}\right)^{\frac{1}{2k-2}}\left( cos\theta + isin\theta \right)$$ 
 
 $$=(b+1)\left( \frac{1}{k^2}\right)^{\frac{k}{2k-2}}cosk\theta + i(b-1)\left( \frac{1}{k^2}\right)^{\frac{k}{2k-2}}sink\theta + ~~~~~~~~~~~~~~~~~~~~~~~~~~~~~~~$$ $$(b+1)\left( \frac{1}{k^2}\right)^{\frac{1}{2k-2}}cos\theta- i(b-1)\left( \frac{1}{k^2}\right)^{\frac{1}{2k-2}}sin\theta$$
 
 $$=(b+1)\left( \frac{1}{k^2}\right)^{\frac{1}{2k-2}}\left[ \left( \frac{1}{k^2}\right)^{\frac{k-11}{2k-2}}cosk\theta + cos\theta \right]+  (b-1)\left( \frac{1}{k^2}\right)^{\frac{1}{2k-2}}\left[ \left( \frac{1}{k^2}\right)^{\frac{k-11}{2k-2}}sink\theta - sin\theta  \right]i $$
 
 $$= (b+1)\left( \frac{1}{k^2}\right)^{\frac{1}{2k-2}}\left( \frac{1}{k}cosk\theta + cos\theta \right)+ (b+1)\left( \frac{1}{k^2}\right)^{\frac{1}{2k-2}}\left( \frac{1}{k}sink\theta - sin\theta \right)i ~~~+ ~~~~~~~~~~~~~~$$ $$ ~~~~~~~~~~~~~~~~~~~~~~~~~~~~~~~~~~~~~~~~~~~~~~ ~~~~~~~~~~~~~~~~~~~~ 2\left( \frac{1}{k^2}\right)^{\frac{1}{2k-2}}\left( \frac{1}{k}sink\theta - sin\theta \right)i $$
Next, let us have the following substitutions: For $$ R= (b+1)\left( \frac{1}{k^2}\right)^{\frac{1}{2k-2}}\left( \frac{k+1}{k}\right) ,~~~  r=(b+1)\left( \frac{1}{k^2}\right)^{\frac{1}{2k-2}}~~~ and ~~~ \phi=k\theta ,$$ we have $$ R-r = (b+1)\left( \frac{1}{k^2}\right)^{\frac{1}{2k-2}}\left(\frac{1}{k}\right),~~~ \frac{R-r}{r}=\frac{1}{k}~~~ and~~~\theta = \frac{R-r}{r}\phi .$$ Therefore, 

$$ Q_{b,b}\left( \Upsilon_{b,b}\right)=(R-r)cos\phi+rcos\left(\frac{R-r}{r}\phi\right) + \left[ (R-r)sin\phi-rsin\left(\frac{R-r}{r}\phi\right)\right]i + $$ $$ \left[ 2\left( \frac{R-r}{b+1}\right)sin\phi-2\left( \frac{r}{b+1} \right) sin\left(  \frac{R-r}{r} \phi \right) \right]i.$$

Hence, the parametric equation of $ Q_{b,b}\left( \Upsilon_{b,b}\right)$ becomes:
$$ \mathcal{X}(\phi)= (R-r)cos\phi+rcos\left(\frac{R-r}{r}\phi\right)~~~~~~~~~~~~~~~~~~~~~~~~~~~~~~~~~~~~~~~~~~~~~~~~~~~~~~~~~$$ 

$$ \mathcal{Y}(\phi) = \left[ (R-r)sin\phi-rsin\left(\frac{R-r}{r}\phi\right)\right] + \left[ 2\left( \frac{R-r}{b+1}\right)sin\phi-2\left( \frac{r}{b+1} \right) sin\left(  \frac{R-r}{r} \phi \right) \right]. $$ This is the parametric equation of hypocycloid centered at $\left[ 2\left( \frac{R-r}{b+1}\right)sin\phi-2\left( \frac{r}{b+1} \right) sin\left(  \frac{R-r}{r} \phi \right) \right].$\\ Moreover, $$ \frac{R}{r}= \frac{(b+1)\left( \frac{1}{k^2}\right)^{\frac{1}{2k-2}}\left( \frac{k+1}{k}\right)}{(b+1)\left( \frac{1}{k^2}\right)^{\frac{1}{2k-2}}}=\frac{k+1}{k}.$$ That means, the ratio $\frac{R}{r}$ does not depend upon  constants $b$ and  $c$, and the entire hypocycloid is traced
for $0 \leq \theta \leq  2\pi.$ Therefore the image of the critical circle under this complex-valued harmonic quadrinomial is a $(k+1,k)$-hypocycloid.

 \begin{figure}[!h]%
    \centering
  {\caption{Image of  $\mathcal{M}_{b,b}e^i\theta$ under  $Q_{b,b}(z)= bz^k+\overline{z}^k+b\overline{z}+z$ for $k=2,k=3,$ and $k=49$ from left to right respectively.}\label{img3}}\includegraphics[width=5cm]{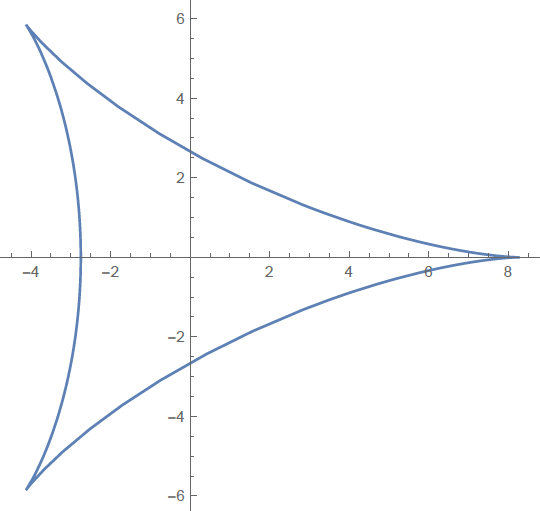}%
   \includegraphics[width=5cm]{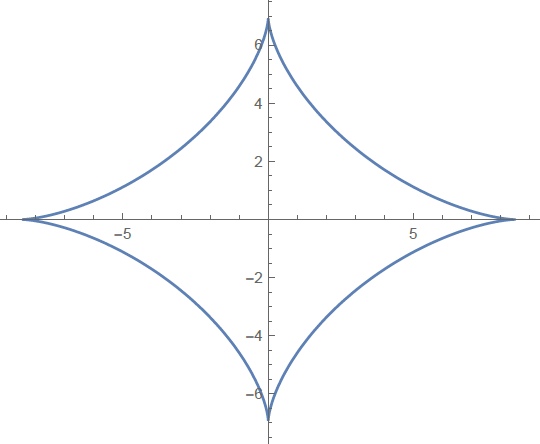}%
    \includegraphics[width=5cm]{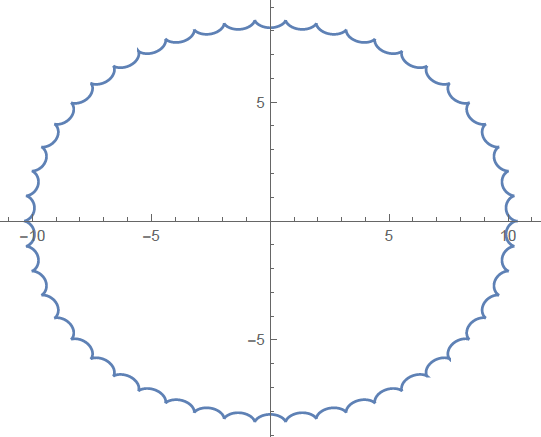}%
\end{figure}
\end{proof}

\begin{theorem} \label{r2}
Let $Q_{b,c}(z)= bz^k+\overline{z}^n+c\overline{z}^m+z$ be a two-parameter family of quadrinomial with $b,c \in \mathbb{R}^+\setminus\{1\},$ $k\geq n>m,$  and $k,n,m \in \mathbb{N}.$  Let $\mathcal{M}_{b,c}:=\left[\frac{(c^2-1)}{k^2(b^2-1)}\right]^{\frac{1}{2k-2}}.$ The quadratic quadrinomial $Q_{b,b}$ has no $\mathcal{M}_{b,b}-$modular roots. 
\end{theorem}
\begin{proof}
An $\mathcal{M}_{b,b}$-modular roots are $\frac{1}{2}$ unit away from the origin. Indeed, $$d(0, \mathcal{M}_{b,b})= \left| \left[\frac{(b^2-1)}{2^2(b^2-1)}\right]^{\frac{1}{2(2)-2}}\right| = \left( \frac{1}{4} \right)^{\frac{1}{2}}= \frac{1}{2}.$$ So, $\mathcal{M}_{b,b}$-modular roots are located on the circle with radius $\frac{1}{2}$ and centered at the origin. 
For $\mathcal{M}_{b,c}$-modular zeros, we have $$ \Upsilon(z)=\left\{z:|z|=\mathcal{M}_{b,c} = \left[\frac{(c^2-1)}{k^2(b^2-1)}\right]^{\frac{1}{2k-2}} \right\} $$ as a critical circle. The corresponding quadratic quadrinomial is $$Q_{b,c}(z)= bz^2+\overline{z}^2+c\overline{z}+z.$$ Our aim is to show that there are no such zeros. Suppose not. Then there is atleast one point $z_0 \in \mathbb{C}$ such that $$bz_0^2+\overline{z_0}^2+b\overline{z_0}+z_0=0.$$ $$\Longrightarrow b \left[ z_0^2+\overline{z_0} \right]= \overline{-z_0}^2-z_0.$$ $$ \Longrightarrow b\left| z_0^2+\overline{z_0} \right| = \left| \overline{z_0}^2+z_0 \right|.$$ $$\Longrightarrow b=\frac{\left| \overline{z_0}^2+z_0 \right|}{\left| z_0^2+\overline{z_0} \right|}.$$ Note that,
$$\left| z_0^2+\overline{z_0} \right| = 0 \Leftrightarrow  z_0^2+\overline{z_0}  = 0 \Leftrightarrow |z_0|^2=|z_0| \Leftrightarrow \frac{1}{4}=\frac{1}{2}.$$ Hence, $\left| z_0^2+\overline{z_0} \right| \neq 0.$ As a result, $$b=\frac{\left| \overline{z_0}^2+z_0 \right|}{\left| z_0^2+\overline{z_0} \right|}=1,$$ which is a contradiction to the fact that $b \neq 1.$ Therefore $\mathcal{M}_{b,b}e^{i\theta}$ is a zero free curve for $k=n=2.$

 \begin{figure}[!h]%
    \centering
   {\caption{Zeros of $Q_{b,b}(z)$ around $\mathcal{M}_{b,b}e^i\theta$ for $b=2$ and $b=12$ respectively.}\label{2}}\includegraphics[width=5cm]{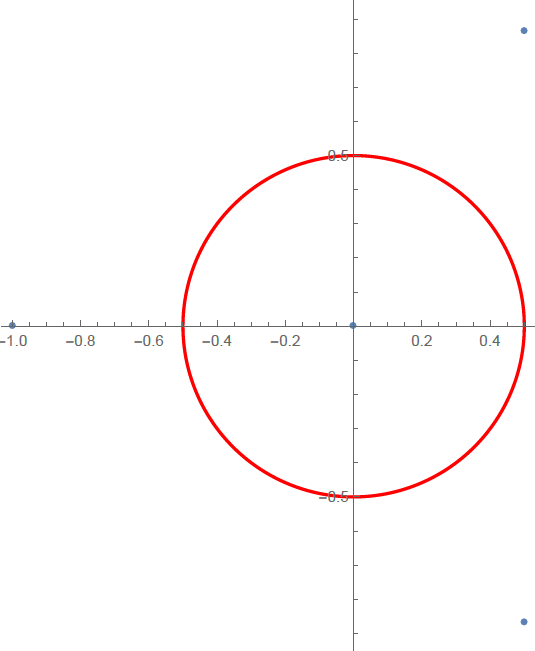}%
   \includegraphics[width=5cm]{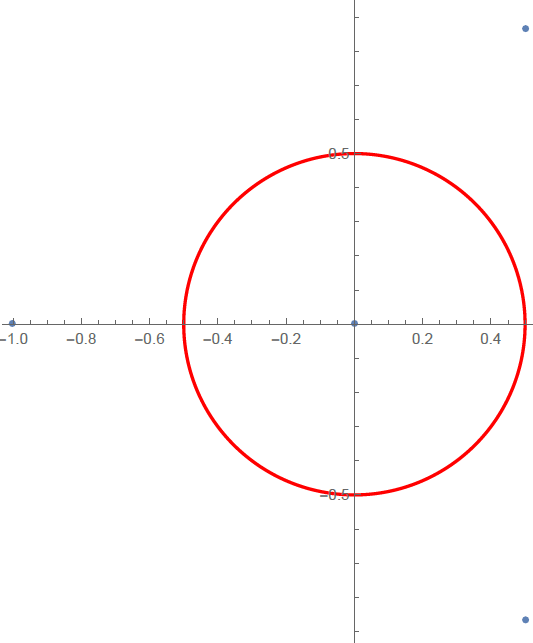}%
\end{figure}
\end{proof}
 
 \section{\textcolor{blue}{Conclusion}}$\label{c}$ 
In this paper, we studied the image of critical circle under complex-valued harmonic quadrinomial and zero free curve for quadratic quadrinomial when the two parameters are equal. We found that the image of a certain critical circles under a complex-valued harmonic quadrinomial of degree $k$ is a hypocycloid with $k+1$ cusps. We also determined that there is no  $\mathcal{M}_{b,b}$-modular zeros. But when the two parameters are different, the type of image as well as whether or not the critical curve contains zeros of quadrinomial remains open.
\section*{Acknowledgments}
We would like to thank Addis Ababa University, Simons Foundation,and ISP(International Science Program) at department level for providing us opportunities and financial support. We also thank Prof. Boris Shapiro and Prof. Rikard B{\"o}gvad for fruitful discussions during our stay at Stockholm University, department of Mathematics.
\section*{Declaration of Interest of Statement}
The authors declare that there are no conflicts of interest regarding the publication of this paper.

\end{document}